\documentclass[12pt, oneside]{article}
\usepackage{latexsym}
\usepackage{graphicx}
\usepackage{setspace}
\usepackage[dvips, ps, rotate, all, color]{xy}
\usepackage{geometry}
\usepackage[psamsfonts]{amsfonts}
\usepackage{amsmath, amsthm, amssymb}
\usepackage{enumerate}
\usepackage[usenames]{color}
\usepackage[sc, center, small]{titlesec}

\xyoption{arc}
\titlelabel{\thetitle. \quad}

\usepackage{verbatim}

\newcommand{\IZ}{\mathbb Z}

\newcommand{\IR}{\mathbb R}
\newcommand{\IC}{\mathbb C}

\newcommand{\hyperH}{\mathbb H}

\newcommand{\lsup}[2]{%
        \ensuremath{{}^{#2}\!{#1}}}

\newtheorem{Thm}{Theorem}
\newtheorem{Cor}[Thm]{Corollary}
\newtheorem{Def}[Thm]{Definition}
\newtheorem{Prop}[Thm]{Proposition}

\newtheorem{Example}[Thm]{Example}
\newtheorem{Remark}[Thm]{Remark}

\advance \topmargin by -\headheight
\advance \topmargin by -\headsep
\evensidemargin \oddsidemargin
\begin{document}


\title{The Zeta Function of a Hypergraph}
\author{Christopher K. Storm\\
  Mathematics Department,\\
  Dartmouth College,\\
  \texttt{cstorm@dartmouth.edu}}
\date{\today}
\maketitle

\begin{center}{\bf \large Abstract}\end{center}
\begin{quote}
We generalize the Ihara-Selberg zeta function to hypergraphs in a natural way.  Hashimoto's factorization results for biregular bipartite graphs apply, leading to exact factorizations.  For $(d,r)$-regular hypergraphs, we show that a modified Riemann hypothesis is true if and only if the hypergraph is Ramanujan in the sense of Winnie Li and Patrick Sol\'e.  Finally, we give an example to show how the generalized zeta function can be applied to graphs to distinguish non-isomorphic graphs with the same Ihara-Selberg zeta function.
\end{quote}

\bigskip

\section{Introduction}

The aim of this paper is to give a non-trivial generalization of the Ihara-Selberg zeta function to hypergraphs and show how our generalization can be thought of as a zeta function on a graph.  We will be concerned with producing generalizations of many of the results known for the Ihara-Selberg zeta function: factorizations, functional equations in specific cases, and an interpretation of a ``Riemann hypothesis."  We will also look at some of the properties of hypergraphs that are determined by our generalization.

Later in this section, we will give the appropriate hypergraph definitions and path definitions necessary for the zeta function.  Keqin Feng and Winnie Li give an Alon-Boppana type result for the eigenvalues of the adjacency operator of hypergraphs \cite{FLi} which will motivate a definition for \emph{Ramanujan hypergraphs} given by Li and Sol\'e \cite{LiS}.  We will also give the appropriate definitions to define a ``prime cycle" in a hypergraph and give a formal definition of the zeta function.

Section 2 is concerned with generalizing a construction of Motoko Kotani and Toshikazu Sunada \cite{Sun}.  The prime cycles in the hypergraph will correspond exactly to \emph{admissible cycles} in a strongly connected, oriented graph.  This will let us write the zeta function as a determinant involving the \emph{Perron-Frobenius operator} $T$ of the strongly connected, oriented graph.  The zeta function will look like $\det(I-uT)^{-1}$, which is a rational function of the form one divided by a polynomial.

In Section 3 we explore in more detail the connection between a hypergraph and its \emph{associated bipartite graph} and what happens as prime cycles are represented in the bipartite graph.  This will allow us to realize the zeta function in terms of the Ihara-Selberg zeta function of the bipartite graph.  Theorem \ref{Thm:HypergraphZetaAsBipartite} details this connection in full. We remark that our generalization is non-trivial in the sense that there are infinitely many hypergraphs whose generalized zeta function is never the Ihara-Selberg zeta function of a graph.  We then get very nice factorization results from Ki-Ichiro Hashimoto's work \cite{H}, found in Theorem \ref{Thm:RegularHypergraphGeneralFactor}.  As corollaries to Hashimoto's factorization results, we will be able to give functional equations and connect the Riemann hypothesis to the Ramanujan condition for a hypergraph.  Theorem \ref{Thm:RamanujanEquivalence} shows that a Riemann hypothesis is true if and only if the hypergraph is Ramanujan.  We will also show how our zeta function fits into hypergraph theory and can give information about whether a hypergraph is \emph{unimodular} and about some coloring properties for the hypergraph.  These results are not new but more a matter of framing previously known work in this context.

Finally, in Section 4 we show how this generalization can actually be applied to graphs.  One impediment to the Ihara-Selberg zeta function being truly useful as a graph invariant is that two $k$-regular graphs are \emph{cospectral}---their adjacency operators have the same spectrum---if and only if they have the same zeta function \cite{Mellein, Q}.  We will examine an example of two $3$-regular graphs constructed by Harold Stark and Audrey Terras \cite{ST2} which have the same zeta function but can be shown explicitly to be non-isomorphic by computing our zeta function in an appropriate way.

For the rest of this section, we fix our terminology and definitions.  For the most part, we are following \cite{FLi, LiS} for our definitions.  A hypergraph $\hyperH = (V, E)$ is a set of \emph{hypervertices} $V$ and a set of \emph{hyperedges} $E$ where each hyperedge is a nonempty set whose elements come from $V$, and the union of all the hyperedges is $V$.  We note that a hypervertex may not be repeated in the same hyperedge; although, with appropriate care it is easy to generalize to this case.  We allow hyperedges to repeat.  A hypervertex $v$ is \emph{incident} to a hyperedge $e$ if $v \in e$.  Finally, we call the cardinality of a hyperedge $e$, denoted $|e|$, the \emph{order} of the hyperedge.

Using the incidence relation, we can associate a bipartite graph $B$ to $\hyperH$ in the following way: the vertices of $B$ are indexed by $V(\hyperH)$ and $E(\hyperH)$.  Vertices $v \in V(\hyperH)$ and $e \in E(\hyperH)$ are \emph{adjacent} in $B$ if $v$ is incident to $e$.  Given a hypergraph $\hyperH$, we will denote by $B_\hyperH$ the bipartite graph formed in this manner.  Given a hypergraph $\hyperH$, we can construct its dual $\hyperH^*$ by letting its hypervertex set be indexed by $E(\hyperH)$ and its hyperedges by $V(\hyperH)$.  We can use the bipartite graph to then construct the appropriate incidence relation.

The associated bipartite graph is a very important tool in the study of hypergraphs.  For now, we can use it to define an \emph{adjacency matrix} for $\hyperH$.  The adjacency matrix $A$ is a matrix whose rows and columns are parameterized by $V(\hyperH)$.  The $ij$-entry of $A$ is the number of directed paths in $B_\hyperH$ from $v_i$ to $v_j$ of length $2$ with no backtracking.  

The adjacency matrix is symmetric---given a path of length 2 from $v_i$ to $v_j$, we traverse it backwards to get a path from $v_j$ to $v_i$---so it has real eigenvalues.  We denote these eigenvalues, referred to as a set as the \emph{spectrum of the adjacency matrix}, by $\lambda_1, \cdots, \lambda_{|V(\hyperH)|}$.  The \emph{spectrum of $\hyperH$} is defined to be the spectrum of $A$ and satisfies
\[ \Delta \geq \lambda_1 \geq \lambda_2 \geq \cdots \geq \lambda_{|V(\hyperH)|} \geq -\Delta\]
for some $\Delta \in \IR$.
\begin{Def}
A hypergraph $\hyperH$ is $(d, r)$-regular if:
\begin{enumerate}
\item Every hypervertex is incident to $d$ hyperedges, and
\item Every hyperedge contains $r$ hypervertices.
\end{enumerate}
\end{Def}

For a $(d, r)$-regular hypergraph, we have $\lambda_1 = d(r-1)$, and the fundamental question becomes how large can the other eigenvalues be?  Feng and Li, generalizing a technique of Alon Nilli \cite{N}, give the following Alon-Boppana type result to address this question \cite{FLi}:
\begin{Thm}[Feng and Li]
Let $\{\hyperH_m\}$ be a family of connected $(d, r)$-regular hypergraphs with $|V(\hyperH_m)| \rightarrow \infty$ as $m \rightarrow \infty$.  Then
\[\liminf \, \lambda_2(\hyperH_m) \geq r - 2 + 2\sqrt{q} \, \, \, {\rm as} \, m \rightarrow \infty,\]
where $q = (d - 1)(r - 1) = d(r-1) - (r - 1)$.
\label{Thm:HyperBoppana}
\end{Thm}
Theorem \ref{Thm:HyperBoppana} is the key ingredient for defining Ramanujan hypergraphs; however, we need to explore the connection between $\hyperH$, $B_\hyperH$, and $\hyperH^*$ a bit more before we give the definition.  When $\hyperH$ is $(d, r)$-regular, we also have that $\hyperH^*$ is $(r, d)$-regular.  Then we can relate the adjacency operators of $\hyperH$, $B_\hyperH$, and $\hyperH^*$ as follows:
\begin{equation}A(B_\hyperH) = \begin{pmatrix} 0&M \\ \lsup{M}{t}&0 \end{pmatrix},
\label{eq:BipartiteAdjacencyOperator} \end{equation}
\begin{equation}A(B_\hyperH)^2 = \begin{pmatrix} M \lsup{M}{t}&0 \\ 0&\lsup{M}{t} M \end{pmatrix} = 
			\begin{pmatrix} A(\hyperH) + dI_V&0 \\ 0&A(\hyperH^*) + rI_E\end{pmatrix},\label{eq:BipartiteAdjacencySquared}\end{equation}
where $M = M(V, E)$ is the incidence matrix of $\hyperH$, and $I_V$ and $I_E$ are identity matrices with rows and columns parameterized by $V$ and $E$ respectively.  Eq. (\ref{eq:BipartiteAdjacencyOperator}) follows from the definitions of the associated bipartite graph $B_\hyperH$ and by ordering the vertices in $B_\hyperH$ in the same way as the hypervertices and hyperedges of $\hyperH$.  To see Eq. (\ref{eq:BipartiteAdjacencySquared}), we first note that the $(i, j)$-entry of $A(B_\hyperH)^k$ is the number of paths of length $k$ from $v_i$ to $v_j$ \cite{West}.  Hence, the $(i, j)$-entry of $A(B_\hyperH)^2$ is the number of paths of length $2$ from $v_i$ to $v_j$ without backtracking plus the number of paths of length $2$ from $v_i$ to $v_j$ with backtracking.  The adjacency operators of $\hyperH$ and $\hyperH^*$ account for the paths without backtracking.  The only way to have a path of length $2$ from $v_i$ to $v_j$ with backtracking is for $i$ and $j$ to be equal.  Then, the number of such paths is either $d$ or $r$, depending on if $v_i$ comes from a hypervertex or a hyperedge, respectively, in $\hyperH$.  This accounts for the identity terms in the expression.

We let $P(x)$, $P^*(x)$, and $Q(x)$ denote the characteristic polynomials of $A(\hyperH)$, $A(\hyperH^*)$, and $A(B_\hyperH)^2$ respectively.  Then by Eq. (\ref{eq:BipartiteAdjacencySquared}), the characteristic polynomials are related by
\begin{equation}
Q(x) = P(x - d)P^*(x - r).
\label{eq:CharRelation}
\end{equation}
Since the eigenvalues of $A(B_\hyperH)^2$ are all non-negative, this relation forces the eigenvalues of $\hyperH$ and $\hyperH^*$ to be at least $-d$ and $-r$ respectively.  We can also relate $P(x)$ and $P^*(x)$ directly as shown in \cite{CDS}:
\begin{equation}
x^{|V|}P^*(x - r) = x^{|E|}P(x - d).
\label{eq:HyperCharRelation}
\end{equation}

This gives a very explicit connection between the spectra of $\hyperH$ and $\hyperH^*$.  When $d$ and $r$ are not equal, comparing the powers of $x$ in both sides of Eq. (\ref{eq:HyperCharRelation}) gives the \emph{obvious} eigenvalue $-d$ of $\hyperH$ with multiplicity $|V(\hyperH)| - |E(\hyperH)|$ or $-r$ of $\hyperH^*$ with multiplicity $|E(\hyperH)| - |V(\hyperH)|$, depending on whether $d < r$ or $r < d$.

Taking into account potential obvious eigenvalues and Theorem \ref{Thm:HyperBoppana}, we define Ramanujan hypergraphs:
\begin{Def}[Li and Sol\'e]
Let $\hyperH$ be a finite, connected $(d, r)$-regular hypergraph.  We say $\hyperH$ is a \emph{Ramanujan hypergraph} if
\begin{equation}
|\lambda - r + 2| \leq 2\sqrt{(d-1)(r-1)},
\end{equation}
for all \emph{non-obvious} eigenvalues $\lambda \in \text{Spec}(\hyperH)$ such that $\lambda \neq d(r-1).$
\label{Def:RamHypergraph}
\end{Def}

This will be the basics of what we need for general hypergraph definitions.  We refer the interested reader to \cite{Ber1, Ber2, FLi, LiS} for more information on hypergraphs and their spectra.  We also point out that there are other potential definitions for Ramanujan hypergraphs that depend on the operators one wishes to study \cite{Li}.  For some explicit constructions of Ramanujan hypergraphs of the type treated here, we refer the reader to \cite{MST}.  We now turn our attention to the definition of the generalized Ihara-Selberg zeta function of a hypergraph.  We recommend the series of articles by Harold Stark and Audrey Terras to the reader interested in current theory on Ihara-type zeta functions on graphs \cite{ST, ST2, ST3}.  Recently, there have also been a number of generalizations of the zeta functions to digraphs as well as buildings \cite{MS1, MS2, D}.

\begin{figure}
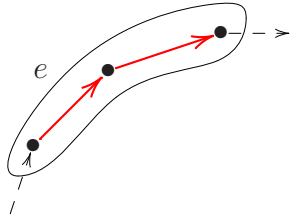

\[ \xy
(0,0)*{\bullet} = "v1";
(10, 10)*{\bullet} = "v2";
(25, 15)*{\bullet} = "v3";
(-3, -9)*{} = "B";
(34, 15)*{} = "C";
(-2, -4)*{} = "A";
(1, 10)*{e};
{\ar@{-->} "B"; "v1"};
{\ar@{-->} "v3"; "C"};
{\ar@[red]@[thicker]@{->} "v1"; "v2"};
{\ar@[red]@[thicker]@{->} "v2"; "v3"};
"A"; "A" **\crv{ (-6, -4) & (5, 15) & (27, 20) & (29, 18) & (27, 10) & (10, 10) & (3, -4) };
\endxy \]
\caption{Hyperedge backtracking in a $3$-edge $e$.}
\label{fig:HyperedgeBacktrack}
\end{figure}

To define our zeta function, we need the appropriate concept of a ``prime cycle."  A \emph{closed path} in $\hyperH$ is a sequence $c = (v_1, e_1, v_2, e_2, \cdots, v_k, e_k, v_1)$, of \emph{length} $k = |c|$, such that $v_i \in e_{i - 1}, e_i$ for $i \in \IZ / k\IZ$.  Note that this implies that $v_1 \in e_k$ so that this path really is ``closed."  We say $c$ has  \emph{hyperedge backtracking} if there is a subsequence of $c$ of the form $(e, v, e)$.  If we have hyperedge backtracking, this means that we use a hyperedge twice in a row.  In general, when we exclude cycles with hyperedge backtracking, it will be permissible to return directly to a hypervertex so long as a different hyperedge is used.  We give an example of hyperedge backtracking in Figure \ref{fig:HyperedgeBacktrack}.  We denote by $c^m$ the $m$-multiple of $c$ formed by going around the closed path $m$ times.  Then, $c$ is \emph{tail-less} if $c^2$ does not have hyperedge backtracking.  If, in addition to having no hyperedge backtracking and being tail-less, $c$ is not the non-trivial $m$-multiple of some other closed path $b$, we say that $c$ is a \emph{primitive cycle}.  Finally, we can impose an equivalence relation on primitive cycles via cyclic permutation of the sequence that defines the cycles.  We call a representative of $[c]$ a \emph{prime cycle}.  We note that direction of travel does matter, so given a triangle in a graph, it can actually be viewed as two prime cycles.

We now define the generalized Ihara-Selberg zeta function of a hypergraph:
\begin{Def}
For $u \in \IC$ with $|u|$ sufficiently small, we define the \emph{generalized Ihara-Selberg zeta function} of a finite hypergraph $\hyperH$ by
\[\zeta_\hyperH(u) = \prod_{\mathfrak{p} \in P}\left(1 - u^{|\mathfrak{p}|}\right)^{-1},\]
where $P$ is the set of prime cycles of $\hyperH$.
\label{Def:ZetaDef}
\end{Def}
\begin{Remark}
A graph $X$ can be viewed as a hypergraph where every hyperedge has order $2$.  In this case, the definitions we've given---and in particular the definition for hyperedge backtracking---correspond exactly to those needed to define prime cycles in graphs.  The zeta function $\zeta_X(u)$ is, then, exactly the Ihara-Selberg zeta function $Z_X(u)$.
\end{Remark}

In the next section, we will focus on giving an initial factorization of $\zeta_\hyperH(u)$, which represents the zeta function as a determinant of explicit operators.  In Section 3, we show more explicit factorizations, using results of Hyman Bass \cite{B} and Hashimoto \cite{H}.  Finally, in Section 4, we give an interpretation of this zeta function as a graph zeta function and show how it can distinguish non-isomorphic graphs that are cospectral.

\begin{center}{\bf Acknowledgments}\end{center}
The author would like to thank Dorothy Wallace and Peter Winkler for several valuable discussions and comments in preparing this manuscript.

\section{The Oriented Line Graph Construction}

The goal of this section is to generalize the construction of an ``oriented line graph" which Kotani and Sunada \cite{Sun} use to begin factoring the Ihara-Selberg zeta function.  The idea is to start with a hypergraph and construct a strongly connected, oriented graph which has the same cycle structure.  By changing the problem from hypergraphs to strongly connected, oriented graphs we will actually make finding an explicit expression for $\zeta_\hyperH(u)$ much simpler.

We first define some terms for oriented graphs.  For an oriented graph, an oriented edge $e = \{x, y\}$ is an ordered pair of vertices $x$, $y \in V$.  We say that $x$ is the \emph{origin} of $e$, denoted by $o(e)$, and $y$ is the \emph{terminus} of $e$, denoted by $t(e)$.  We also have the \emph{inverse edge} $\bar{e}$ given by switching the origin and terminus.  We say that an oriented, finite graph $X^o = (V, E^o)$ is \emph{strongly connected} if, for any $x, y \in V$, there exists an \emph{admissible path} $c$ with $o(c) = x$ and $t(c) = y$.  A path $c = (e_1, \cdots, e_k)$ is \emph{admissible} if $e_i \in E^o$ and $o(e_i) = t(e_{i-1})$ for all $i$.  We say that $o(c) = o(e_1)$ and $t(c) = t(e_k)$.

Let $\hyperH$ be a finite, connected hypergraph.  We label the edges of $\hyperH$: $E =\\ \{e_1, e_2, \cdots, e_m\}$ and fix $m$ colors $\{c_1, c_2, \cdots, c_m\}$.  We now construct an edge-colored graph $G\hyperH_c$ as follows.  The vertex set $V(G\hyperH_c)$ is the set of hypervertices $V(\hyperH)$.  For each hyperedge $e_j \in E(\hyperH)$, we construct a $|e_j|$-clique in $G\hyperH_c$ on the hypervertices in $e_j$ by adding an edge, joining $v$ and $w$, for each pair of hypervertices $v, w \in e_j$.  We then color this $|e_j|$-clique $c_j$.  Thus if $e_j$ is a hyperedge of order $i$, we have $i \choose 2$ edges in $G\hyperH_c$, all colored $c_j$.

Once we've constructed $G\hyperH_c$, we arbitrarily orient all of the edges.  We then include the inverse edges as well, so we finish with a graph $G\hyperH_c^o$ which has twice as many colored, oriented edges as $G\hyperH_c$.

\begin{figure}
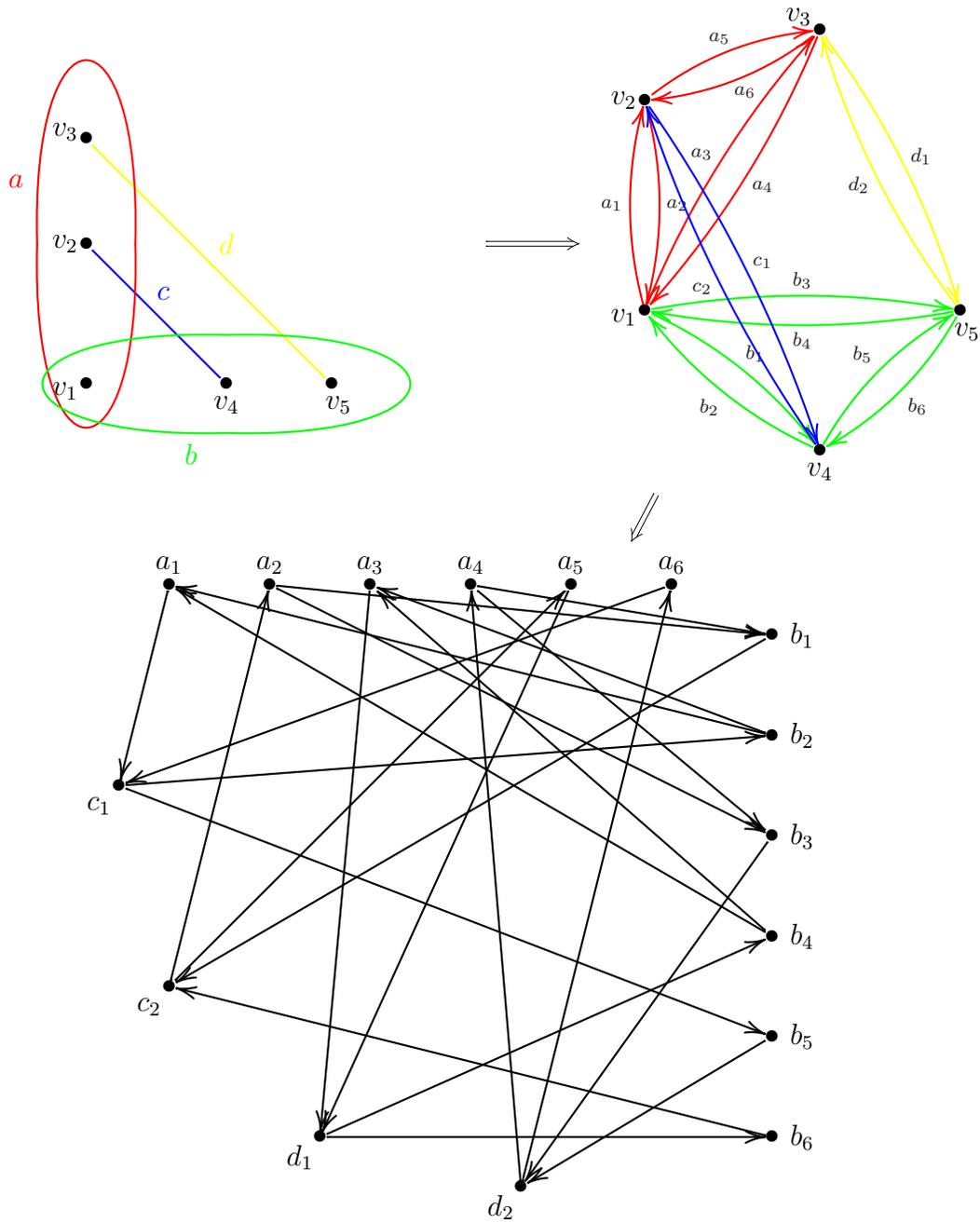

\[ \xy
(-5,0)*!C++[thicker]\xybox{
(0,0)*{\bullet} = "v1" +(-3, -1)*{v_1};
(0,20)*{\bullet} = "v2" +(-3, 0)*{v_2};
(0,35)*{\bullet} = "v3" +(-3, 1)*{v_3};
(20,0)*{\bullet} = "v4" +(0, -3)*{v_4};
(35,0)*{\bullet} = "v5" +(1, -3)*{v_5};
(-10, 29)*[red]{a};
(15, -10)*[green]{b};
(11, 13)*[blue]{c};
(20, 20)*[yellow]{d};
(-7, 20)*{} = "A";
(7, 20)*{} = "B";
(20, -7)*{} = "C";
(20, 7)*{} = "D";
{\ar@{-}@[blue] "v2"; "v4"};
{\ar@{-}@[yellow] "v3"; "v5"};
"A"; "B" **[red]\crv{ (-8, 55) & (8, 55) };
"A"; "B" **[red]\crv{ (-8, -15) & (8, -15) };
"C"; "D" **[green]\crv{ (-15, -8) & (-15, 8)};
"C"; "D" **[green]\crv{ (55, -8) & (55, 8)};
} = "x";
(75,0)*!C++[thicker]\xybox{
(0,0)*{\bullet} = "v1" +(-3, -1)*{v_1};
(0,30)*{\bullet} = "v2" +(-3, 0)*{v_2};
(25,40)*{\bullet} = "v3" +(-3, 2)*{v_3};
(25,-20)*{\bullet} = "v4" +(0, -3)*{v_4};
(45,0)*{\bullet} = "v5" +(1, -3)*{v_5};
{\ar^{a_1}@{->}@/^/@[red] "v1"; "v2"};
{\ar^{a_2}@{->}@/^/@[red] "v2"; "v1"};
{\ar^{a_3}@{->}@/^/@[red] "v1"; "v3"};
{\ar^{a_4}@{->}@/^/@[red] "v3"; "v1"};
{\ar^{a_5}@{->}@/^/@[red] "v2"; "v3"};
{\ar^{a_6}@{->}@/^/@[red] "v3"; "v2"};
{\ar^{b_1}@{->}@/^/@[green] "v1"; "v4"};
{\ar^{b_2}@{->}@/^/@[green] "v4"; "v1"};
{\ar^{b_3}@{->}@/^/@[green] "v1"; "v5"};
{\ar^{b_4}@{->}@/^/@[green] "v5"; "v1"};
{\ar^{b_5}@{->}@/^/@[green] "v4"; "v5"};
{\ar^{b_6}@{->}@/^/@[green] "v5"; "v4"};
{\ar^{c_1}@{->}@/^/@[blue] "v2"; "v4"};
{\ar^{c_2}@{->}@/^/@[blue] "v4"; "v2"};
{\ar^{d_1}@{->}@/^/@[yellow] "v3"; "v5"};
{\ar^{d_2}@{->}@/^/@[yellow] "v5"; "v3"};
} = "y";
{\ar@{=>} "x"; "y"};
(27,-92)*!C++[black][thicker]\xybox{0;/r.34pc/:
(0,0)*{\bullet}="a1"+(0, 2)*{a_1};
(10,0)*{\bullet}="a2"+(0, 2)*{a_2};
(20,0)*{\bullet}="a3"+(0, 2)*{a_3};
(30,0)*{\bullet}="a4"+(0, 2)*{a_4};
(40,0)*{\bullet}="a5"+(0, 2)*{a_5};
(50,0)*{\bullet}="a6"+(0, 2)*{a_6};
(60,-5)*{\bullet}="b1"+(3, 0)*{b_1};
(60,-15)*{\bullet}="b2"+(3, 0)*{b_2};
(60,-25)*{\bullet}="b3"+(3, 0)*{b_3};
(60,-35)*{\bullet}="b4"+(3, 0)*{b_4};
(60,-45)*{\bullet}="b5"+(3, 0)*{b_5};
(60,-55)*{\bullet}="b6"+(3, 0)*{b_6};
(-5,-20)*{\bullet}="c1"+(-2, -2)*{c_1};
(0,-40)*{\bullet}="c2"+(-2, -2)*{c_2};
(15,-55)*{\bullet}="d1"+(-2, -2)*{d_1};
(35,-60)*{\bullet}="d2"+(-2, -2)*{d_2};
{\ar@{->} "a1"; "c1"};
{\ar@{->} "a2"; "b3"};
{\ar@{->} "a2"; "b1"};
{\ar@{->} "a3"; "d1"};
{\ar@{->} "a4"; "b3"};
{\ar@{->} "a4"; "b1"};
{\ar@{->} "a5"; "d1"};
{\ar@{->} "a6"; "c1"};
{\ar@{->} "b1"; "c2"};
{\ar@{->} "b2"; "a3"};
{\ar@{->} "b2"; "a1"};
{\ar@{->} "b3"; "d2"};
{\ar@{->} "b4"; "a3"};
{\ar@{->} "b4"; "a1"};
{\ar@{->} "b5"; "d2"};
{\ar@{->} "b6"; "c2"};
{\ar@{->} "c1"; "b5"};
{\ar@{->} "c1"; "b2"};
{\ar@{->} "c2"; "a2"};
{\ar@{->} "c2"; "a5"};
{\ar@{->} "d1"; "b4"};
{\ar@{->} "d1"; "b6"};
{\ar@{->} "d2"; "a4"};
{\ar@{->} "d2"; "a6"};
} = "z";
{\ar@{=>} "y"; "z"};
\endxy \]
\caption{We begin with a hypergraph $\hyperH$, already colored, in the top left.  Then we construct one possible edge-colored oriented graph $G\hyperH_c^o$.  From this graph, we construct the corresponding oriented line graph.  We notice that there are no edges that go from $a_i$ to $a_j$; this is because they represent the red edges in $G\hyperH_c^o$.}
\label{fig:HypergraphOrientedLineSmaller}
\end{figure}

Finally, we construct the \emph{oriented line graph} $\hyperH_L^o = (V_L, E_L^o)$ associated with our choice of $G\hyperH_c^o$ by
\begin{align*}
V_L &= E(G\hyperH_c^o), \\
E_L^o &= \{(e_i, e_j) \in E(G\hyperH_c^o) \times E(G\hyperH_c^o); c(e_i) \neq c(e_j), t(e_i) = o(e_j)\},
\end{align*}
where $c(e_i)$ is the colored assigned to the oriented edge $e_i \in E(G\hyperH_c^o)$.  If our hypergraph $\hyperH$ was a graph to begin with, for any oriented edge $e \in E(G\hyperH_c^o)$, the only oriented edge with the same color is $\bar{e}$.  Then, the oriented line graph construction given here is exactly that given by Kotani and Sunada \cite{Sun}.  See Figure \ref{fig:HypergraphOrientedLineSmaller} for an example of this construction.
\begin{Prop}
Suppose $\hyperH$ is a finite, connected hypergraph where each hypervertex is in at least two hyperedges and which has more than two prime cycles.  Then, the oriented line graph $H_L^o$ is finite and strongly connected.
\end{Prop}
\begin{proof}
The vertices of $H_L^o$ are of the form $\{v, w\}_{e}$ where $e \in E(\hyperH)$ and $v, w \in e$.  This catalogues using the hyperedge $e$ to go from $v$ to $w$.  To show that $H_L^o$ is strongly connected, we must show that given two subsequences $\{v_1, e_1, v_2\}$ and $\{v_k, e_k, v_{k+1}\}$ with $e_1, e_k \in E(\hyperH)$, $v_1, v_2 \in e_1$, and $v_k, v_{k+1} \in e_k$, there exists a path $c$ in $\hyperH$ of the form $c = (v_1, e_1, v_2, e_2, \cdots, e_{k-1}, v_k, e_k, v_{k+1})$ such that $c$ has no hyperedge backtracking.  Since $c$ has no hyperedge backtracking, we can use this path to construct a path in $H_L^o$ which starts at $\{v_1, v_2\}_{e_1}$ and finishes at $\{v_k, v_{k+1}\}_{e_k}$.

Since $\hyperH$ is connected and every hypervertex is in at least $2$ hyperedges, there exists a path with no hyperedge backtracking $d$ which begins with $(v_1, e_1, v_2, \cdots)$ and finishes at vertex $v_k$.  Now there are two cases.  Either the path used $e_k$ in the last step to get to $v_k$ or it did not.  If the path did not use $e_k$, we can use $e_k$ to go to $v_{k+1}$, and we are done.  In the second case, we need the additional hypothesis that there are more than two prime cycles.  We can get the desired path by leaving $v_k$ via a hyperedge different than $e_k$.  Then there is some cycle (which may have a tail) which returns to $v_k$ via the other hyperedge.  Then we can go from $v_k$ to $v_{k+1}$ via $e_k$.  This yields the desired path.  In essence, we need more than two prime cycles to allow ourselves to ``turn around" if we get going in the wrong direction.  Hence, $H_L^o$ is strongly connected.

That $H_L^o$ is finite is clear since $\hyperH$ is finite.
\end{proof}

For $m \geq 1 \in \IZ$, we let $N_m$ be the number of admissible closed paths of length $m$ in $\hyperH_L^o$.  Then, we define the \emph{zeta function} of $\hyperH_L^o$ by
\begin{equation}
Z^o_{\hyperH_L^o}(u) = \exp\left(\sum_{m = 1}^{\infty}\frac{1}{m}N_mu^m \right).
\label{Def:Zeta}
\end{equation}

The initial factorization for this zeta function is determined in terms of the \emph{Perron-Frobenius operator} $T: C(V_L) \mapsto C(V_L)$ given by
\[(Tf)(x) = \sum_{e \in E_0(x)}f(t(e)),\]
where $E_0(x) = \{e \in E_o | o(e) = x\}$ is the set of all oriented edges with $x$ as their origin vertex.

Kotani and Sunada \cite{Sun} give the details to let us factor $Z^o_{\hyperH_L^o}(u)$ in terms of its Perron-Frobenius operator:
\begin{Thm}[Kotani and Sunada]
Suppose $\hyperH_L^o$ is a finite, oriented graph which is strongly connected and not just a circuit.  Then
\[Z^o_{\hyperH_L^o}(u) = \exp\left(\sum_{m = 1}^{\infty}\frac{1}{m}N_mu^m \right) = \det(I - uT)^{-1},\]
where $T$ is the Perron-Frobenius operator of $\hyperH_L^o$.
\label{thm:OrientedZetaFactor}
\end{Thm}
\begin{proof}
We only sketch the details:
\begin{enumerate}
\item Convergence in a disk about the origin follows from the Perron-Frobenius theorem \cite{Gant}.
\item The factorization was essentially given by Rufus Bowen and O. E. Lanford III \cite{Bow}.
\end{enumerate}
\end{proof}

We denote by $P_L^o$ the set of admissible prime cycles; then, we have the following Euler Product expansion
\[Z^o_{\hyperH_L^o}(u) = \prod_{\mathfrak{p} \in P_L^o}(1 - u^{|\mathfrak{p}|})^{-1}\]
which is Theorem 2.3 in \cite{Sun}.  Viewing the zeta function in this manner, we need only show a correspondence between the prime cycles of $\hyperH$ and the admissible prime cycles of $\hyperH_L^o$:
\begin{Prop}
There is a one-to-one correspondence between prime cycles of length $l$ in $\hyperH$ and admissible prime cycles of length $l$ in $\hyperH_L^o$.  In particular, the zeta function of $\hyperH$ can be written as
\[\zeta_\hyperH(u) = \det(I - uT)^{-1},\]
where $T$ is the Perron-Frobenius operator of $\hyperH_L^o$.
\end{Prop}
\begin{proof}
We show the stated cycle correspondence; then, the factorization will follow from the Euler Product expansion of $Z^o_{\hyperH_L^o}(u)$ and Theorem \ref{thm:OrientedZetaFactor}.

To show the cycle correspondence, we will actually show that there is a correspondence between paths in $\hyperH$ with no hyperedge backtracking and admissible paths in $\hyperH_L^o$.  The cycle correspondence will then follow since all the relations imposed on paths are the same.

Suppose $v$ and $w$ are hypervertices contained in a hyperedge $e$.  Then we denote by $\{v, w\}_e$ the oriented edge in $G\hyperH_c^o$ with origin $v$, terminus $w$, and color given by the color chosen for $e$.  We let $c = (v_1, e_1, v_2, e_2, \cdots, v_k, e_k, v_{k+1})$ be a path in $\hyperH$ with no hypervertex backtracking.  This corresponds to the path $c^o = (\{v_1, v_2\}_{e_1}, \{v_2, v_3\}_{e_2}, \cdots, \{v_k, v_{k+1}\}_{e_k})$ in $G\hyperH_c^o$.  Since there is no hyperedge backtracking, i.e. $e_i \neq e_{i+1}$ at every step, we change colors as we follow each oriented edge.  Then the corresponding path $\tilde{c} = ((\{v_1, v_2\}_{e_1}, \{v_2, v_3\}_{e_2}), (\{v_2, v_3\}_{e_2}, \{v_3, v_4\}_{e_3}), \cdots, \\(\{v_{k-1}, v_k\}_{e_{k-1}}, \{v_k, v_{k+1}\}_{e_k})) $ in $\hyperH_L^o$ is admissible with length $k$.

Similarly, given an admissible path in $\hyperH_L^o$, we can realize it as a path in $G\hyperH_c^o$ which changes colors at every step.  That means the corresponding path in $\hyperH$ changes hyperedges at every step; i.e., that it does not have hyperedge backtracking.  Then lengths, then, are the same.
\end{proof}

In particular, this theorem means that the zeta function is a rational function and provides a tool to make some initial calculations.  To get more precise factorizations, we shall look more closely at the relationship between a hypergraph and its associated bipartite graph.

\section{Further Factorizations}

In the last section, we were able to realize the generalized Ihara-Selberg zeta function as a determinant of explicit operators.  In this section, we will see that by shifting our view to the associated bipartite graph of a hypergraph, we can do much better.  Once we've established the relation between cycles in  hypergraphs and cycles in bipartite graphs that we need, we will draw very heavily from Hashimoto's work on zeta functions of bipartite graphs \cite{H}.  To help keep clear what structure we are referring to, we will continue to call cycles in a hypergraph cycles but will call cycles in the associated bipartite graph geodesics.

To motivate the relation we are looking for, we look at a simple example.  In Figure \ref{fig:PathInHyperToBipartite}, we look at the primitive cycle given by $c = (v_1, e_1, v_2, e_3, v_4, e_2, v_1)$.  This corresponds to a primitive geodesic $\tilde{c} = (v_1, \{v_1, e_1\}, e_1, \{e_1, v_2\}, \cdots, \{e_2, v_1\}, v_1)$ in the associated bipartite graph.  In fact, this sort of correspondence is true in general:
\begin{figure}
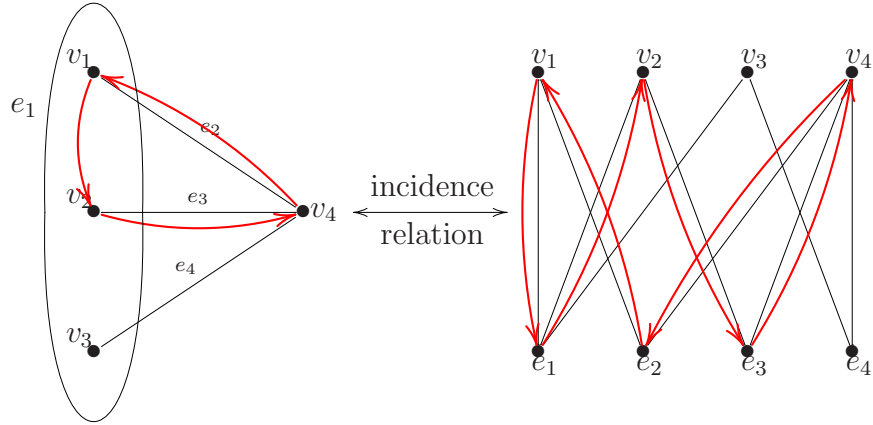

\[\xy0;/r.22pc/:
(0,0)*++{\xy
(0,20)*{\bullet}="v1"+(-2, 2)*{v_1};
(0,0)*{\bullet}="v2"+(-2, 2)*{v_2};
(0,-20)*{\bullet}="v3"+(-2, 2)*{v_3};
(30,0)*{\bullet}="v4"+(3, 0)*{v_4};
(-7, 0)*{} = "A";
(7, 0)*{} = "B";
(-10, 15)*{e_1};
{\ar^{e_2}@{-} "v1"; "v4"};
{\ar^{e_3}@{-} "v2"; "v4"};
{\ar^{e_4}@{-} "v3"; "v4"};
{\ar@/_/@[red]@[thicker]@{->} "v1"; "v2"};
{\ar@/_/@[red]@[thicker]@{->} "v2"; "v4"};
{\ar@/_/@[red]@[thicker]@{->} "v4"; "v1"};
"A"; "B" **\crv{ (-8, 40) & (8, 40) };
"A"; "B" **\crv{ (-8, -40) & (8, -40) };
\endxy} = "x";
(75, 0)*++{\xy
(0,20)*{\bullet}="v1"+(1, 2)*{v_1};
(15,20)*{\bullet}="v2"+(1, 2)*{v_2};
(30,20)*{\bullet}="v3"+(1, 2)*{v_3};
(45,20)*{\bullet}="v4"+(1, 2)*{v_4};
(0,-20)*{\bullet}="e1"+(1, -2)*{e_1};
(15,-20)*{\bullet}="e2"+(1, -2)*{e_2};
(30,-20)*{\bullet}="e3"+(1, -2)*{e_3};
(45,-20)*{\bullet}="e4"+(1, -2)*{e_4};
{\ar@{-} "v1"; "e1"};
{\ar@{-} "v1"; "e2"};
{\ar@{-} "v2"; "e1"};
{\ar@{-} "v2"; "e3"};
{\ar@{-} "v3"; "e1"};
{\ar@{-} "v3"; "e4"};
{\ar@{-} "v4"; "e2"};
{\ar@{-} "v4"; "e3"};
{\ar@{-} "v4"; "e4"};
{\ar@/_/@[red]@[thicker]@{->} "v1"; "e1"};
{\ar@/_/@[red]@[thicker]@{->} "e1"; "v2"};
{\ar@/_/@[red]@[thicker]@{->} "v2"; "e3"};
{\ar@/_/@[red]@[thicker]@{->} "e3"; "v4"};
{\ar@/_/@[red]@[thicker]@{->} "v4"; "e2"};
{\ar@/_/@[red]@[thicker]@{->} "e2"; "v1"};
\endxy}= "y";
{\ar@{<->} "x"; "y"};
(37, 4)*{\text{incidence}};
(37, -3)*{\text{relation}};
\endxy
\]
\caption{An example of a primitive cycle of length $3$ in a hypergraph and a corresponding primitive geodesic of length $6$ in its associated bipartite graph.}
\label{fig:PathInHyperToBipartite}
\end{figure}
\begin{Prop}
Let $\hyperH$ be a finite, connected hypergraph with associated bipartite graph $B_\hyperH$.  Then there is a one-to-one correspondence between prime cycles of length $l$ in $\hyperH$ and prime geodesics of length $2l$ in $B_\hyperH$.
\end{Prop}
\begin{proof}
We will begin with a representative of a prime cycle of length $l$ in $\hyperH$.  Let $c = (v_1, e_1, \cdots, v_l, e_l, v_1)$ be a primitive cycle in $\hyperH$.  Then we claim that $\\ \tilde{c} = (v_1, \{v_1, e_1\}, e_1, \cdots, v_l, \{v_l, e_l\}, e_l, \{e_l, v_1\}, v_1)$ is a primitive geodesic in $B_\hyperH$.  It is clear that $\tilde{c}$ is both closed and primitive if $c$ is, so we need only check to be sure $\tilde{c}$ has no backtracking or tails.

Let's look at what hyperedge backtracking in the hypergraph means.  We say that $c$ has hyperedge backtracking if we use the same hyperedge twice in a row.  If we think about the bipartite graph side, this means we leave a vertex in the set from $E(\hyperH)$, go to a vertex in the set $V(\hyperH)$ and then backtrack to the vertex in $E(\hyperH)$.  Still on the bipartite side, the only other way to backtrack is to go from a vertex in $V(\hyperH)$ to a vertex in $E(\hyperH)$ and directly back.  Thus, we would have the following sequence in the hypergraph: $(v_i, e_i, v_i)$.  This type of sequence is expressly disallowed unless $v_i$ is repeated more than once in $e_i$.  If this happens, there is a multiple edge in $B_\hyperH$ representing this, which means we can actually return to the first vertex without backtracking.  Putting all of this together, we see that not hyperedge backtracking in $\hyperH$ is equivalent to not backtracking on the corresponding path in $B_\hyperH$.  Once we know that backtracking isn't an issue, having no tails also follows immediately since backtracking in $\tilde{c}^2$ would correspond to hyperedge backtracking in $c^2$.  Thus, each prime cycle of length $l$ in $\hyperH$ corresponds to a prime geodesic of length $2l$ in $B_\hyperH$.

We now look at prime geodesics in $B_\hyperH$ and show that they correspond to prime cycles in $\hyperH$.  Without loss of generality, we can assume that the first entry in a representative of a prime geodesic in $B_\hyperH$ is a vertex parameterized by the set $V(\hyperH)$.  If it is not, we simply shift the cycle one slot in either direction, and we will have an appropriate representative because the graph is bipartite.  Suppose the representative looks like $\tilde{c} = (v_1, \{v_1, e_1\}, e_1, \cdots, v_l, \{v_l, e_l\}, e_l, \{e_l, v_1\}, v_1)$; then we have the following primitive cycle in $\hyperH$: $c = (v_1, e_1, \cdots, v_l, e_l, v_1)$.  This is a primitive cycle by the same reasons as above since $\tilde{c}$ is a primitive geodesic.  Also, $|\tilde{c}| = 2l = 2|c|$, so we see that given a prime geodesic in $B_\hyperH$, it corresponds to a prime cycle of half the length in $\hyperH$.
\end{proof}

This correspondence means that we can relate the generalized Ihara-Selberg zeta function of a hypergraph to the Ihara-Selberg zeta function of its associated bipartite graph.
\begin{Thm}
Let $\hyperH$ be a finite, connected hypergraph such that every hypervertex is in at least two hyperedges.  Then,
\[ \zeta_\hyperH(u) = Z_{B_\hyperH}(\sqrt{u}).\]
\label{Thm:HypergraphZetaAsBipartite}
\end{Thm}
\begin{proof}
Let $P_\hyperH$ be the set of prime cycles on $\hyperH$ and $P_{B_\hyperH}$ the set of prime geodesics on $B_\hyperH$.  Then we rely on the previous proposition to write:
\begin{align*}
\zeta_\hyperH(u) &= \prod\limits_{\mathfrak{p} \in P_\hyperH}\left(1 - u^{|\mathfrak{p}|}\right)^{-1}
			   = \prod\limits_{\mathfrak{p} \in P_\hyperH}\left(1 - {u^{2|\mathfrak{p}|}}^{\frac{1}{2}}\right)^{-1} \\
			   &= \prod\limits_{\ell \in P_{B_\hyperH}}\left(1 - {u^{|\ell|}}^{\frac{1}{2}}\right)^{-1}
			   = Z_B(\sqrt{u}).
\end{align*}
\end{proof}

As an immediate consequence of this relation, we see that, for an arbitrary hypergraph $\hyperH$ which satisfies the conditions of Theorem \ref{Thm:HypergraphZetaAsBipartite} we can relate its zeta function to the zeta function of its dual hypergraph $\hyperH^*$.
\begin{Cor}
Suppose $\hyperH$ satisfies the conditions of Theorem \ref{Thm:HypergraphZetaAsBipartite}.  Then,
\[\zeta_\hyperH(u) = \zeta_{\hyperH^*}(u)\]
\end{Cor}
\begin{proof}
$\hyperH$ and $\hyperH^*$ have the same associated bipartite graph, by definition.  Then we apply Theorem \ref{Thm:HypergraphZetaAsBipartite}.
\end{proof}

In addition, we can rewrite Hyman Bass's Theorem \cite{B} on factoring the zeta function of a graph to give us a form of $\zeta_\hyperH(u)$ which is more amenable to computation.  We first state Bass's Theorem:
\begin{Thm}[Bass]
Let $X$ be a finite, connected graph with adjacency operator $A$ and operator $Q$ defined by $D - I$ where $D$ is the diagonal operator with the degree of vertex $v_i$ in the $i$th slot of the diagonal.  Let $I$ be the $|V| \times |V|$ identity operator.  Then,
\[Z_{X}(u) = (1 - u^2)^{\chi(X)} \det(I - uA + u^2Q)^{-1}\]
where $\chi = |V| - |E|$ is the Euler Number of the graph $X$.
\label{Thm:Bass}
\end{Thm}

Given a hypergraph $\hyperH$, we apply Theorem \ref{Thm:Bass} to factor $Z_\hyperH(u)$, giving us a computable factorization of $\zeta_\hyperH(u)$:
\begin{Cor}
Let $\hyperH$ be a finite, connected hypergraph such that every hypervertex is in at least two hyperedges.  Let $A_{B_\hyperH}$ be the adjacency operator on $B_\hyperH$, and let $Q_{B_\hyperH}$ be the operator on $B_\hyperH$ defined by $D - I$ where $D$ is the diagonal operator with the degree of vertex $v_i$ in the $i$th slot of the diagonal.  Let $I$ be the $m \times m$ identity operator where $m = |V(\hyperH)| + |E(\hyperH)|$.  Then
\[\zeta_\hyperH(u) = Z_{B_\hyperH}(\sqrt{u}) = (1 - u)^{\chi(B_\hyperH)} \det(I - \sqrt{u}A_{B_\hyperH} + uQ_{B_\hyperH})^{-1},\]
where $\chi(B_\hyperH) = |V(B_\hyperH)| - |E(B_\hyperH)|$.
\label{Cor:BassForGeneralizedIS}
\end{Cor}
\begin{Remark}
We make a few notes, highlighting how we can compute each of the terms that show up in Corollary \ref{Cor:BassForGeneralizedIS}:
\begin{enumerate}
\item Despite the square root that appears in this factorization, $\zeta_\hyperH(u)$ is a rational function.  We see this clearly in the previous section, but we can also recover it quickly by recalling that a bipartite graph only has prime cycles of even length.
\item The adjacency operator of $B_\hyperH$ can be quickly constructed from the incidence matrix of $\hyperH$ as in Eq. (\ref{eq:BipartiteAdjacencyOperator}).
\item Similarly, we can construct the operator $Q_{B_\hyperH}$ quickly by considering the degrees of vertices in the associated bipartite graph.  If $x$ is a vertex which comes from $V(\hyperH)$, we have $d(x)$ is the number of hyperedges of which $x$ is a member, counting possible multiplicity.  If $x$ comes from $E(\hyperH)$, then $d(x)$ is the order of the associated hyperedge.  From these two facts, we can easily reconstruct $Q_{B_\hyperH}$.
\item We also see that $|V(B_\hyperH)| = |V(\hyperH)| + |E(\hyperH)|$.  In addition, $|E(B_\hyperH)|$ can be directly computed in two different ways via
\[|E(B_\hyperH)| = \sum_{e \in E(\hyperH)} |e| = \sum_{v \in V(\hyperH)} \sharp\{e \in E(\hyperH); v \in e\}.\]
\end{enumerate}
\end{Remark}
\begin{Example}
We compute the generalized Ihara-Selberg zeta function of the hypergraph in Figure \ref{fig:PathInHyperToBipartite} in two ways.  By going through the oriented line graph, we have
\[\zeta_\hyperH(u)^{-1} =  \det(I - uT) = (1-u)(1 + u + u^2 - 5u^3 - 5u^4 - 5u^5 + 4u^6 + 4u^7 + 4u^8).\]
We can also compute the zeta function of the associated bipartite graph by using Bass's Theorem  to realize
\[Z_{B_\hyperH}(u)^{-1} = (1-u^2)(1 + u^2 + u^4 - 5u^6 - 5u^8 - 5u^{10} + 4u^{12} + 4u^{14} + 4u^{16}).\]
Then we can directly see the relation $\zeta_\hyperH(u) = Z_{B_\hyperH}(\sqrt{u})$.
\label{exa:ComputeZeta}
\end{Example}

We emphasize that Corollary \ref{Cor:BassForGeneralizedIS} is mainly useful for computation.  In general, the diagonal entries of the $Q$ matrix will not all be the same, making it quite difficult to manipulate the factorization for theoretical results.  Theorem \ref{Thm:HypergraphZetaAsBipartite} makes it clear that the problem of factoring the generalized zeta function is really a problem of factoring the zeta function of a bipartite graph.  Fortunately, in \cite{H}, Hashimoto deals with this question in great detail.

We reformulate Hashimoto's Main Theorem(III) \cite{H} into our context to get the following theorem:
\begin{Thm}
Suppose that $\hyperH$ is a finite, connected $(d, r)$-regular hypergraph with $d \geq r$.  Let $n_1 = |V(\hyperH)|$, $n_2 = |E(\hyperH)|$, and $q = (d-1)(r-1)$.  Let $A$ be the adjacency operator of $\hyperH$, and let $A^*$ be the adjacency operator of $\hyperH^*$.  Then one has
\begin{align*}
\zeta_\hyperH&(u)^{-1} \\
&= (1-u)^{-\chi(B_\hyperH)}(1+(r-1)u)^{(n_2 - n_1)}\times\det[I_{n_1} - (A - r + 2)u +qu^2] \\
&= (1-u)^{-\chi(B_\hyperH)}(1+(d-1)u)^{(n_1 - n_2)}\times\det[I_{n_2} - (A^* - d +2)u +qu^2],
\end{align*}
where $-\chi(B_\hyperH) = n_1(d-1) - n_2 = n_2(r-1) - n_1$.
\label{Thm:RegularHypergraphGeneralFactor}
\end{Thm}

Theorem \ref{Thm:RegularHypergraphGeneralFactor} will provide the tool we need to produce theoretical results about the generalized Ihara-Selberg zeta function on $(d, r)$-regular hypergraphs.  The condition that $d \geq r$ is actually not a problem.  If $\hyperH$ is a $(d, r)$-regular hypergraph; then, $\hyperH^*$ is $(r, d)$-regular.  Thus, if $d < r$, we simply consider $\hyperH^*$ as our starting point instead.

\subsection{Consequences of the Factorization}

Our first observation is that the zeta function of a hypergraph is a non-trivial generalization of the Ihara-Selberg zeta function.  By this, we mean that we can produce an infinite number of zeta functions which are not the zeta function of any graph.  A simple way to produce zeta functions which did not come from a graph is encoded in the next proposition.
\begin{Prop}
Suppose $X$ is a finite graph with no vertices of degree $1$, and $\hyperH$ is a finite hypergraph with every hypervertex in at least $2$ hyperedges.  Then,
\begin{enumerate}
\item The degree of the polynomial $Z_X(u)^{-1}$ is $2|E(X)|$.
\item The degree of the polynomial $\zeta_\hyperH(u)^{-1}$ is $\sum\limits_{e \in E(\hyperH)} |e|$.
\end{enumerate}
\label{Prop:Degrees}
\end{Prop}
\begin{proof}
\begin{enumerate}
\item Let $X$ be a finite graph with no vertices of degree $1$.  Then by Bass's Theorem \cite{B},
\[Z_X(u)^{-1} = (1-u^2)^{|E| - |V|} \times \det(I - uA + u^2Q).\]

The degree of the determinant term is $2|V|$, and the degree of the explicit polynomial is $2|E(X)| - 2|V(X)|$.  Hence, the degree of $Z_X(u)^{-1}$ is $2|E(X)| - 2|V(X)| + 2|V(X)| = 2|E(X)|$.
\item Let $\hyperH$ be a finite hypergraph with associated bipartite graph $B_\hyperH$.  Then by Theorem \ref{Thm:HypergraphZetaAsBipartite},
\[\zeta_\hyperH(u)^{-1} = Z_{B_\hyperH}(\sqrt{u})^{-1}.\]

From the previous part, we see that the degree of $Z_{B_\hyperH}(\sqrt{u})^{-1}$ is $|E(B_\hyperH)|$.  We can compute this explicitly as $|E(B_\hyperH)| = \sum\limits_{e \in E(\hyperH)}|e|$.
\end{enumerate}
\end{proof}

If a graph $X$ has a vertex of degree $1$, we can remove that vertex and the edge to which it is adjacent without changing the zeta function.  By removing all of these types of vertices until we are left with a graph with every vertex having degree at least $2$, we see that the zeta function of the graph we started with will be the zeta function of a graph that satisfies Proposition \ref{Prop:Degrees}.  Hence, the inverse of the zeta function of a graph will always have even degree.  If we wish to exhibit hypergraphs with zeta functions that did not arise from some graph, we need only find a hypergraph for which $\sum\limits_{e \in E(\hyperH)}|e|$ is odd.
\begin{Example}
In Example \ref{exa:ComputeZeta}, we computed the zeta function of the hypergraph appearing in Figure \ref{fig:PathInHyperToBipartite}.  We see that the inverse of the zeta function has odd degree, so this is an example of a hypergraph which produces a zeta function that no graph could produce.
\end{Example}

Before we turn to a discussion of the poles of the zeta function of a $(d, r)$-regular hypergraph, we look at some of the symmetry that Hashimoto's factorization gives us.  These functional equations are in the spirit of those given by Stark and Terras in \cite{ST}.
\begin{Cor}
Suppose that $\hyperH$ is a finite connected $(d, r)$-regular hypergraph with $d \geq r$.  Let $n_1 = |V(\hyperH)|$, $n_2 = |E(\hyperH)|$, $q = (d-1)(r-1)$, and $\chi = \chi(B_\hyperH)$.  Let $A$ be the adjacency operator of $\hyperH$, and let $A^*$ be the adjacency operator of $\hyperH^*$.  Finally, suppose $p(u)$ is a polynomial in $u$ that satisfies $p(u)^\eta = \pm(qu^2)^\eta p(\frac{1}{qu})^\eta$.  Then we have the following functional equations for $\zeta_\hyperH(u)$:
\begin{enumerate}
\item $\Lambda_\hyperH(u) = p(u)^{n_1}(1-u)^{-\chi}(1+(r-1)u)^{n_2-n_1}\zeta_\hyperH(u) = \pm\Lambda_\hyperH(\frac{1}{qu})$.
\item $\tilde{\Lambda}_\hyperH(u) = p(u)^{n_2}(1-u)^{-\chi}(1+(d-1)u)^{n_1-n_2}\zeta_\hyperH(u) = \pm\tilde{\Lambda}_\hyperH(\frac{1}{qu})$.
\end{enumerate}
\label{Cor:HypergraphFunctionalEquations}
\end{Cor}
\begin{proof}
The strategy is really one of brute force factorization, using Theorem \ref{Thm:RegularHypergraphGeneralFactor}.  By Theorem \ref{Thm:RegularHypergraphGeneralFactor}, we can write $\zeta_\hyperH(u)$ as
\[\zeta_\hyperH(u) = (1-u)^\chi(1+(r-1)u)^{(n_1 - n_2)}\times\det[I_{n_1} - (A - r + 2)u +qu^2]^{-1}.\]
Substituting this expression into $\Lambda_\hyperH(u)$, we have
\[\Lambda_\hyperH(u) = p(u)^{n_1}\times\det[I_{n_1} - (A - r + 2)u +qu^2]^{-1}.\]
We now algebraically manipulate the determinant term:
\begin{align*}
\det[I_{n_1} - (A - r + 2)u +qu^2]^{-1} &= \det[\frac{qu^2}{qu^2} - (A - r + 2)\frac{qu^2}{qu} + \frac{qu^2}{1}]^{-1} \\
&= \left(\frac{1}{qu^2}\right)^{n_1} \times \det[\frac{1}{qu^2} - (A - r +2)\frac{1}{qu} + \frac{1}{1}]^{-1} \\
&= \left(\frac{1}{qu^2}\right)^{n_1} \times \det[\frac{1}{1}I_{n_1} - (A - r +2)\frac{1}{qu} + \frac{q}{(qu)^2}]^{-1}.
\end{align*}
We substitute this back into the expression for $\Lambda_\hyperH(u)$ and then use the given condition for $p(u)^{n_1}$:
\begin{align*}
\Lambda_\hyperH(u) &= p(u)^{n_1} \times \left(\frac{1}{qu^2}\right)^{n_1} \times \det[\frac{1}{1}I_{n_1} - (A - r +2)\frac{1}{qu} + \frac{q}{(qu)^2}]^{-1} \\
&= \pm(qu^2)^{n_1}p(\frac{1}{qu})^{n_1} \times \left(\frac{1}{qu^2}\right)^{n_1} \times \det[\frac{1}{1}I_{n_1} - (A - r +2)\frac{1}{qu} + \frac{q}{(qu)^2}]^{-1} \\
&= \pm p(\frac{1}{qu})^{n_1}\times\det[\frac{1}{1}I_{n_1} - (A - r +2)\frac{1}{qu} + \frac{q}{(qu)^2}]^{-1} \\
&= \pm\Lambda_\hyperH(\frac{1}{qu}).
\end{align*}
This completes the first functional equation.  The second one is identical, using Hashimoto's second factorization.  We leave it as an exercise to the reader.
\end{proof}
\begin{Remark}
Using Corollary \ref{Cor:HypergraphFunctionalEquations}, we can write down several explicit functional equations for $(d, r)$-hypergraphs with $d \geq r$.
\begin{enumerate}
\item $\Lambda_\hyperH(u) = (1-u)^{n_1-\chi}(1+(r-1)u)^{(n_2-n_1)}(1-qu)^{n_1}\zeta_\hyperH(u) = \Lambda_\hyperH(\frac{1}{qu})$.
\item $\tilde{\Lambda}_\hyperH(u) = (1-u)^{n_2-\chi}(1+(d-1)u)^{(n_1-n_2)}(1-qu)^{n_2}\zeta_\hyperH(u) = \tilde{\Lambda}_\hyperH(\frac{1}{qu})$.
\item $\Xi_\hyperH(u) = (1-u)^{-\chi}(1+(r-1)u)^{(n_2-n_1)}(1+qu^2)^{n_1}\zeta_\hyperH(u) = \Xi_\hyperH(\frac{1}{qu}).$
\item $\tilde{\Xi}_\hyperH(u) = (1-u)^{-\chi}(1+(d-1)u)^{(n_1-n_2)}(1+qu^2)^{n_2}\zeta_\hyperH(u) = \tilde{\Xi}_\hyperH(\frac{1}{qu}).$
\end{enumerate}
\end{Remark}

Now that we have several established functional equations, we turn to the next important question for a zeta function.  We will look at the location of the poles and show that they very explicitly detect the Ramanujan condition on a $(d, r)$-regular hypergraph.  

We assume throughout that $\hyperH$ is a $(d, r)$-regular hypergraph with $d \geq r$.  We let $n_2 = |E(\hyperH)|$, $n_1 = |V(\hyperH)|$, $q = (d-1)(r-1)$, and $A$ be the adjacency operator on $\hyperH$.  Then, we have that $n_2 \geq n_1$ since $d \geq r$.  By Eq. (\ref{eq:HyperCharRelation}), $\hyperH$ has no \emph{obvious eigenvalues} $-d$.  This will simplify our consideration of the Ramanujan condition on $\hyperH$.

We now want to focus on the determinant term in Hashimoto's factorization.  Since $A$ is symmetric, it is diagonalizable, so suppose $Q$ diagonalizes $A$.  Then,
\begin{align*}
\det[I_{n_1} - (A - r& + 2)u +qu^2] = \det\left(Q[I_{n_1} - (A - (r - 2)I_{n_1})u +qu^2I_{n_1}]Q^{-1}\right) \\
&= \det[QI_{n_1}Q^{-1} - (QAQ^{-1} - (r-2)QI_{n_1}Q^{-1})u + qu^2QI_{n_1}Q^{-1}] \\
&= \det[I_{n_1} - (QAQ^{-1} - r+2)u + qu^2] \\
&= \prod_{\lambda \in \text{Spec}(\hyperH)}[1 - (\lambda - r + 2)u + qu^2].
\end{align*}

This is the factorization we need to fully examine the relation between poles of $\zeta_\hyperH(u)$ and eigenvalues of $\hyperH$.  The next two propositions detail the connection fully.
\begin{Prop}
Suppose $\hyperH$ is a $(d, r)$-regular hypergraph with $d \geq r$.  Then,
\begin{enumerate}
\item $\zeta_\hyperH(u)$ has a pole at $u = 1$ with multiplicity $n_1(d-1) - n_2 = n_2(r-2) - n_1 = -\chi(B_\hyperH)$.
\item $\zeta_\hyperH(u)$ has a pole at $u = -\frac{1}{r-1}$ with multiplicity $n_2 - n_1$.
\end{enumerate}
\end{Prop}
\begin{proof}
The first set of poles is contributed by the factor $(1-u)^{\chi(B_\hyperH)}$ given in Theorem \ref{Thm:RegularHypergraphGeneralFactor}.  The second set is from the factor $(1+(r-1)u)^{(n_1-n_2)}$.
\end{proof}
\begin{Prop}
Suppose $\hyperH$ is a $(d, r)$-regular hypergraph with $d \geq r$.  Let $q = (d-1)(r-1)$, then $\hyperH$ is a Ramanujan hypergraph if and only if the poles of $\det[I_{n_1} - (A - r + 2)u +qu^2]^{-1}$ are distributed as below:
\begin{enumerate}
\item There is a simple pole at $u = 1$ and at $u = \frac{1}{q}$.
\item All other poles lie on the circle in the complex plane given by $|r| = \frac{1}{\sqrt{q}}$.
\end{enumerate}
\end{Prop}
\begin{proof}
Since $\hyperH$ is a $(d, r)$-regular hypergraph, there is an eigenvalue $\lambda = d(r-1)$.  We first rewrite the polynomial for this eigenvalue as 
\begin{align*}
f(u) &= qu^2 - (\lambda - r + 2)u + 1 \\
&= qu^2 - (q+1)u + 1 \\
&= (1-u)(1-qu).
\end{align*}
We can then see the roots at $u = 1$ and at $u = \frac{1}{q}$ as claimed in part 1.  We note that if $\hyperH$ is Ramanujan, the second eigenvalue is strictly smaller than $d(r-1)$, so these poles are simple as claimed.

We now look at the eigenvalues which satisfy $\lambda \neq d(r-1)$.  Then the polynomial $f(u) = qu^2 - (\lambda - r + 2)u + 1$ has roots at
\[u = \frac{(r-2-\lambda) \pm \sqrt{(\lambda - r + 2)^2 - 4q}}{2q}.\]
Then $u$ has $\text{Im}(u) \neq 0$ if and only if $(\lambda - r +2)^2 \leq 4q$.  This is true if and only if $|\lambda - r + 2| \leq 2\sqrt{q}$, which is true if and only if $\hyperH$ is Ramanujan, by Definition \ref{Def:RamHypergraph} (there are no obvious eigenvalues to consider by our assumption on $d$ and $r$).  In this case, we can calculate the modulus of the roots by
\begin{align*}
|u|^2 &= \frac{(\lambda-r+2)^2}{4q^2} + \frac{4q - (\lambda - r + 2)^2}{4q^2} \\
&= \frac{4q}{4q^2} = \frac{1}{q}.
\end{align*}
\end{proof}

This gives us a complete characterization of the relation between the poles of the generalized Ihara-Selberg zeta function and the Ramanujan condition on a hypergraph.  We can rewrite the previous two propositions into a modified Riemann hypothesis.
\begin{Def}
Let $\hyperH$ be a $(d, r)$-regular hypergraph with $d \geq r$ and $q = (d-1)(r-1)$.  We then consider $\zeta_\hyperH(q^{-s})$.  We say that $\zeta_\hyperH(q^{-s})$ satisfies the \emph{modified hypergraph Riemann hypothesis} if and only if for
\[\text{Re }s \in (0, 1), \, \, \, \frac{(1+(r-1)q^{-s})^{(n_2 - n_1)}}{\zeta_\hyperH(q^{-s})} = 0 \Longrightarrow \text{Re } s = \frac{1}{2}.\]
\end{Def}

Then the previous two propositions can be summarized by
\begin{Thm}
For a $(d, r)$-regular hypergraph $\hyperH$, $\zeta_\hyperH(q^{-s})$ satisfies the modified hypergraph Riemann hypothesis if and only if $\hyperH$ is a Ramanujan hypergraph.
 \label{Thm:RamanujanEquivalence}
 \end{Thm}

\subsection{Some Hypergraph Properties}

Before we move on and show how the zeta function can be interpreted as a graph zeta function with a restricted cycle set, we show how some well-known hypergraph properties fit into this framework.  In particular, we will be interested in the case when the zeta function is an even function.  A graph is bipartite if and only if its Ihara-Selberg zeta function is even, and we will see that the generalized zeta function indicates some of the generalizations of ``bipartite" to hypergraphs.  The hypergraph theorems we refer to are all from Chapter 20, Section 3 of Berge \cite{Ber1}.

We let $\hyperH = (V, E)$ be a finite hypergraph.  An \emph{equitable q-colouring} of $\hyperH$ is a partition $(S_1, \cdots, S_q)$ of the hypervertices into $q$ classes such that for each $i \in I$ and for $j, j^{'} \leq q$,
\[-1 \leq |E_i \cap S_j | - |E_i \cap S_{j^{'}}| \leq 1.\]
The smallest number $q \geq 2$ for which there exists an equitable $q$-colouring is the \emph{equitable chromatic number} $\kappa(\hyperH)$ of $\hyperH$.  $\hyperH$ is \emph{unimodular} if for each $S \subset V$, the subhypergraph $\hyperH_S$ admits an equitable bicolouring.  A graph is unimodular if and only if it is bipartite, so this definition is a generalization of bipartite for hypergraphs.

We now look at what it means for the generalized Ihara-Selberg zeta function to be an even function:
\begin{Prop}
Let $\hyperH$ be a hypergraph. Then, $\zeta_\hyperH(u) = \zeta_\hyperH(-u)$ for all $u \in \IC$ if and only if every primitive cycle in $\hyperH$ has even length.
\end{Prop}
\begin{proof}
We consider the power series expansion of the zeta function given in Definition \ref{Def:ZetaDef}.  Then $u$ appears to an odd power if and only if there is a prime cycle of odd length.  Hence, the zeta function must be even on a disk about the origin.  Since it continues to the inverse of a polynomial, it must be even throughout the complex plane.
\end{proof}

This is all we need to reframe several of the results cited in \cite{Ber1}:
\begin{Thm}
Suppose $\hyperH$ is a hypergraph with $\zeta_\hyperH(u)$ an even function.  Then, $\hyperH$ is unimodular.
\end{Thm}
\begin{proof}
This is Theorem 10 in Chapter 20 of \cite{Ber1}.
\end{proof}

\begin{Cor}
Suppose $\hyperH$ is a hypergraph.  Then $\zeta_\hyperH(u)$ is even if and only if each hypergraph $\hyperH^{'}$, defined by taking hyperedges to be subsets of hyperedges of $\hyperH$ and hypervertex set to be the union of all the new hyperedges, satisfies $\kappa(\hyperH^{'}) \leq 2$.
\end{Cor}

We now return to the graph case and see how restricting the set of prime cycles in a graph can give information which is more specific to the graph structure and less dependent on the spectrum of the graph.

\section{The Generalized Zeta Function as Applied to Graphs}

\begin{figure}
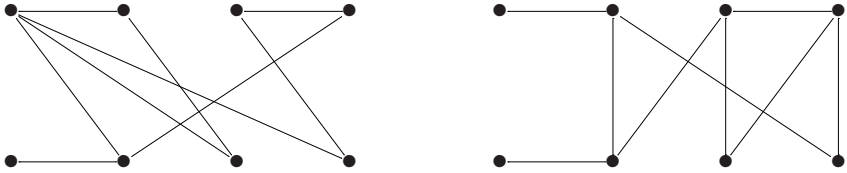

\[ \xy
(0,0)*++{\xy
(0,0)*{\bullet}="v5"; (15,0)*{\bullet}="v6"; (30,0)*{\bullet}="v7"; (45,0)*{\bullet}="v8"; (0,20)*{\bullet}="v1"; (15,20)*{\bullet}="v2"; (30,20)*{\bullet}="v3"; (45,20)*{\bullet}="v4";
{\ar@{-} "v1"; "v2"};
{\ar@{-} "v1"; "v6"};
{\ar@{-} "v1"; "v7"};
{\ar@{-} "v1"; "v8"};
{\ar@{-} "v2"; "v7"};
{\ar@{-} "v3"; "v4"};
{\ar@{-} "v3"; "v8"};
{\ar@{-} "v4"; "v6"};
{\ar@{-} "v5"; "v6"};
\endxy} = "x";
(65, 0)*++{\xy
(0,0)*{\bullet}="v5"; (15,0)*{\bullet}="v6"; (30,0)*{\bullet}="v7"; (45,0)*{\bullet}="v8"; (0,20)*{\bullet}="v1"; (15,20)*{\bullet}="v2"; (30,20)*{\bullet}="v3"; (45,20)*{\bullet}="v4";
{\ar@{-} "v1"; "v2"};
{\ar@{-} "v2"; "v6"};
{\ar@{-} "v2"; "v8"};
{\ar@{-} "v3"; "v4"};
{\ar@{-} "v3"; "v6"};
{\ar@{-} "v3"; "v7"};
{\ar@{-} "v4"; "v8"};
{\ar@{-} "v4"; "v7"};
{\ar@{-} "v5"; "v6"};
\endxy} = "y";
\endxy \]
\caption{Two cospectral graphs with different Ihara-Selberg zeta functions.}
\label{Fig:CospectralDifferentZeta}
\end{figure}

The ideas in this section are motivated, in part, by the question of determining if two given graphs are isomorphic.  We say two graphs are \emph{cospectral} if the spectra of their adjacency matrices are the same.  For general graphs, the Ihara-Selberg zeta function can be useful as a tool for distinguishing graphs since it's possible to have cospectral graphs with different zeta functions.  Figure \ref{Fig:CospectralDifferentZeta} gives an example of cospectral graphs from \cite{HS} which have different zeta functions (there are also examples of graphs which have the same Laplacian spectrum but different zeta functions).  For $k$-regular graphs, however, being cospectral is equivalent to having the same zeta function \cite{Mellein}.  The problem is clearly illustrated by a result of Gregory Quenell \cite{Q}.  We set up some notation before we state his result.

The universal cover of a $k$-regular graph $G$ is the infinite $k$-regular tree, which we denote $X_k$.  We let Aut$(X_k)$ be the group of automorphisms of $X_k$.  Then the graph $G$ can be viewed as the quotient of $X_k$ by a subgroup $H$ of Aut$(X_k)$ that acts freely.  We write $G = H \setminus X_k$; then the vertices of $G$ are the orbits $Hx$ of vertices in $X_k$, and $Hx$ is adjacent to $Hy$ if and only if each element of $Hx$ is adjacent to some element of $Hy$ in $X_k$.  With this framework in mind, we state Quenell's theorem:

\begin{Thm}[Quenell]
Let $H \setminus X_k$ be an $N$-vertex, $k$-regular graph with no loops or parallel edges.  For each integer $n \geq 1$, let
\[P_n = \sum_{[h_i]_H \subset [t_n]_{{\rm Aut}(X_k)}} L(C_H(h_i))\]
where $[t_n]_{{\rm Aut}(X_k)}$ is the Aut$(X_k)$-conjugacy class containing all length-$n$ translations in $H$ and $L(C_H(h_i))$ denotes the length of a generator of the centralizer $C_H(h_i)$ of $h_i$ in $H$.

Then the spectrum of $H\setminus X_k$ determines and is determined by the sequence $P_1, P_2, \cdots, P_N.$
\end{Thm}
We interpret this theorem as saying that the spectrum of the adjacency operator is determined and determines the number of primitive cycles of lengths $1, 2, \cdots, N = |V(X)|$.  These numbers figure prominently in the logarithmic derivative of $Z_X(u)$, giving the zeta function connection.  Hence, if we wish to try to use a zeta function to distinguish cospectral regular graphs, we will need to restrict our paths in some way to try to make them more accurately mimic the unique structure of a given graph.

\begin{figure}
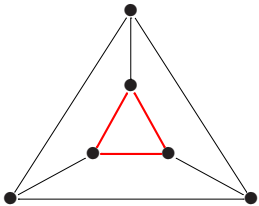

\[\xy
(0,5)*{\bullet} = "v1";
(5, -4)*{\bullet} = "v2";
(-5, -4)*{\bullet} = "v3";
(0, 15)*{\bullet} = "v4";
(16, -10)*{\bullet} = "v5";
(-16, -10)*{\bullet} = "v6";
{\ar@{-} "v1"; "v4"};
{\ar@{-} "v3"; "v6"};
{\ar@{-} "v2"; "v5"};
{\ar@{-} "v4"; "v5"};
{\ar@{-} "v5"; "v6"};
{\ar@{-} "v6"; "v4"};
{\ar@{-}@[red]@[thicker] "v1"; "v2"};
{\ar@{-}@[red]@[thicker] "v2"; "v3"};
{\ar@{-}@[red]@[thicker] "v3"; "v1"};
\endxy\]
\caption{A graph with a triangle singled out.}
\label{fig:GraphTriangle}
\end{figure}
We refer to Figure \ref{fig:GraphTriangle} to illustrate how the generalized Ihara-Selberg zeta function might be used to do this.  We will start with the set of all prime cycles in this graph.  Then we can throw out any prime cycle that uses two red edges in a row.  We actually will be throwing out infinitely many prime cycles when we do this.  We could now define a new zeta function using this smaller set of prime cycles in the same way as before.  It turns out that this is exactly the zeta function for the hypergraph formed by replacing the red triangle with a $3$-edge on the same vertices.  

We could perform the same sort of construction for other graphs by replacing cliques of any size with a hyperedge on the respective vertices.  In this way, we would hope that the path structure would more accurately mirror the structure of the graph and not be as influenced by its spectrum.  For the rest of this section we focus on an example to illustrate this.

\begin{figure}
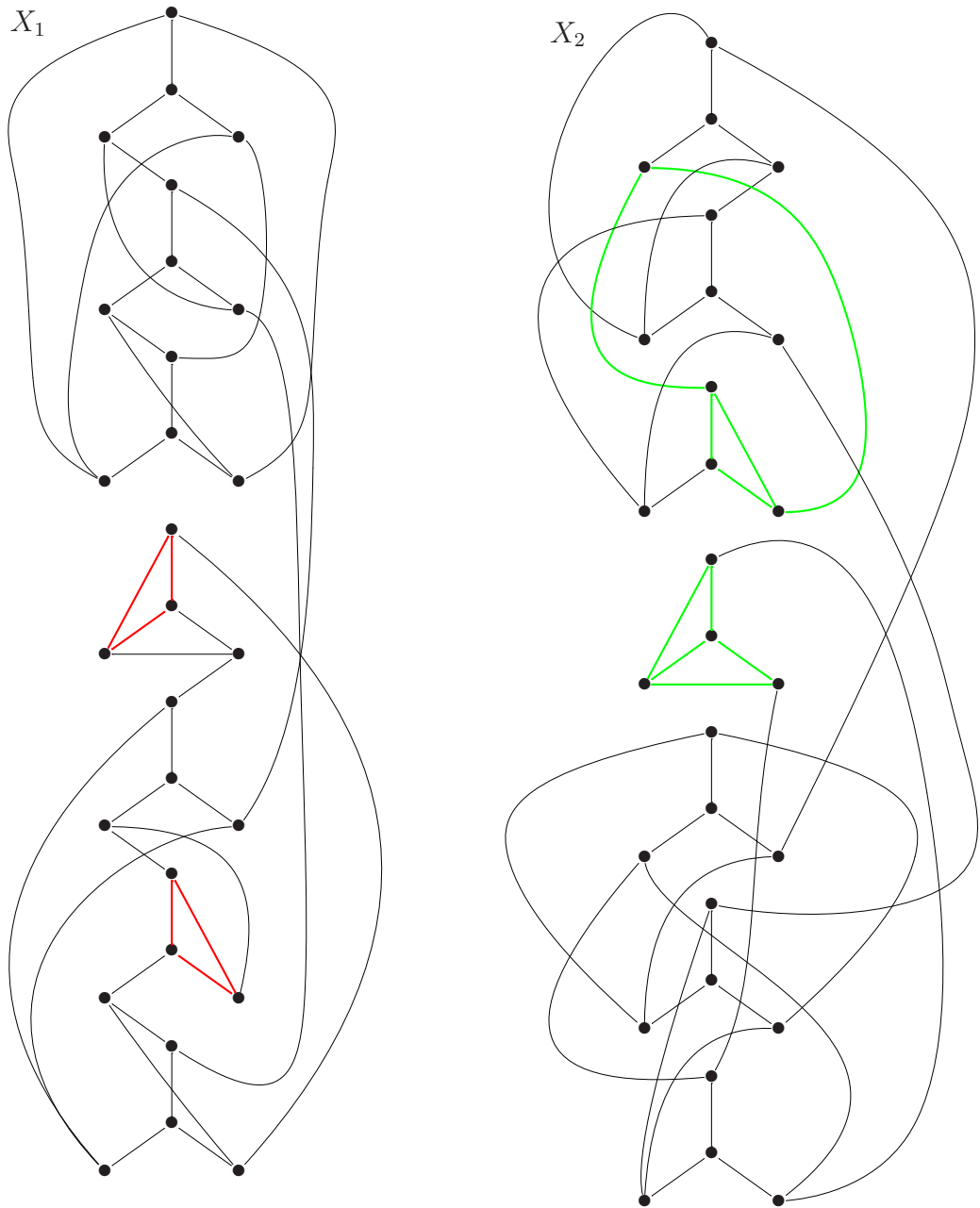

\[\xy0;/r.32pc/:
(0,0)*!C++\xybox{
(0,0)*{\bullet} = "v11";
(-7, -5)*{\bullet} = "v21";
(7, -5)*{\bullet} = "v31";
(0, 8)*{\bullet} = "v41";
(-7, 13)*{\bullet} = "v22";
(7, 13)*{\bullet} = "v32";
(0, 18)*{\bullet} = "v12";
(0, 26)*{\bullet} = "v42";
(0, 36)*{\bullet} = "v13";
(-7, 31)*{\bullet} = "v23";
(7, 31)*{\bullet} = "v33";
(0, 44)*{\bullet} = "v43";
(0, 54)*{\bullet} = "v14";
(-7, 49)*{\bullet} = "v24";
(7, 49)*{\bullet} = "v34";
(0, 62)*{\bullet} = "v44";
(0, 72)*{\bullet} = "v15";
(-7, 67)*{\bullet} = "v25";
(7, 67)*{\bullet} = "v35";
(0, 80)*{\bullet} = "v45";
(0, 90)*{\bullet} = "v16";
(-7, 85)*{\bullet} = "v26";
(7, 85)*{\bullet} = "v36";
(0, 98)*{\bullet} = "v46";
(0, 108)*{\bullet} = "v17";
(-7, 103)*{\bullet} = "v27";
(7, 103)*{\bullet} = "v37";
(0, 116)*{\bullet} = "v47";
(-15, 115)*{X_1};
{\ar@{-} "v11"; "v21"};
{\ar@{-} "v11"; "v31"};
{\ar@{-} "v11"; "v41"};
{\ar@{-} "v12"; "v22"};
{\ar@{-}@[red]@[thicker] "v12"; "v32"};
{\ar@{-}@[red]@[thicker] "v12"; "v42"};
{\ar@{-} "v13"; "v23"};
{\ar@{-} "v13"; "v33"};
{\ar@{-} "v13"; "v43"};
{\ar@{-}@[red]@[thicker] "v14"; "v24"};
{\ar@{-} "v14"; "v34"};
{\ar@{-}@[red]@[thicker] "v14"; "v44"};
{\ar@{-} "v15"; "v25"};
{\ar@{-} "v15"; "v35"};
{\ar@{-} "v15"; "v45"};
{\ar@{-} "v16"; "v26"};
{\ar@{-} "v16"; "v36"};
{\ar@{-} "v16"; "v46"};
{\ar@{-} "v17"; "v27"};
{\ar@{-} "v17"; "v37"};
{\ar@{-} "v17"; "v47"};
{\ar@{-} "v41"; "v22"};
{\ar@{-} "v42"; "v23"};
{\ar@{-} "v43"; "v34"};
{\ar@{-} "v45"; "v26"};
{\ar@{-} "v46"; "v27"};
{\ar@{-}@[red]@[thicker] "v44"; "v24"};
{\ar@{-} "v24"; "v34"};
{\ar@{-}@[red]@[thicker] "v42"; "v32"};
"v47"; "v25" **\crv{ (-20, 110) & (-15, 100) & (-15, 75) & (-13, 70) };
"v47"; "v35" **\crv{ (20, 110) & (15, 100) & (15, 75) & (13, 70) };
"v27"; "v36" **\crv{ (-8, 90) & (0, 85) }; 
"v37"; "v45" **\crv{ (10, 103) & (11, 78) & (5, 80)};
"v37"; "v25" **\crv{ (3, 104) & (-7, 100) & (-10, 85) & (-12, 70) };
"v46"; "v33" **\crv{ (15, 90) & (15, 72) & (14, 40) };
"v26"; "v35" **\crv{ (0, 74)};
"v36"; "v41" **\crv{ (13, 85) & (13, 49) & (15, 10) & (12, 0) };
"v44"; "v31" **\crv{ ( 40, 31) };
"v43"; "v21" **\crv{ (-30, 20) };
"v23"; "v32" **\crv{ (7, 31) & (10, 25) };
"v33"; "v21" **\crv{ (0, 31) & (-18, 20) & (-15, 3) };
"v22"; "v31" **\crv{ (0, 3)};
};
(57,0)*!C++\xybox{
(0,0)*{\bullet} = "v11";
(-7, -5)*{\bullet} = "v21";
(7, -5)*{\bullet} = "v31";
(0, 8)*{\bullet} = "v41";
(-7, 13)*{\bullet} = "v22";
(7, 13)*{\bullet} = "v32";
(0, 18)*{\bullet} = "v12";
(0, 26)*{\bullet} = "v42";
(0, 36)*{\bullet} = "v13";
(-7, 31)*{\bullet} = "v23";
(7, 31)*{\bullet} = "v33";
(0, 44)*{\bullet} = "v43";
(0, 54)*{\bullet} = "v14";
(-7, 49)*{\bullet} = "v24";
(7, 49)*{\bullet} = "v34";
(0, 62)*{\bullet} = "v44";
(0, 72)*{\bullet} = "v15";
(-7, 67)*{\bullet} = "v25";
(7, 67)*{\bullet} = "v35";
(0, 80)*{\bullet} = "v45";
(0, 90)*{\bullet} = "v16";
(-7, 85)*{\bullet} = "v26";
(7, 85)*{\bullet} = "v36";
(0, 98)*{\bullet} = "v46";
(0, 108)*{\bullet} = "v17";
(-7, 103)*{\bullet} = "v27";
(7, 103)*{\bullet} = "v37";
(0, 116)*{\bullet} = "v47";
(-15, 117)*{X_2};
{\ar@{-} "v11"; "v21"};
{\ar@{-} "v11"; "v31"};
{\ar@{-} "v11"; "v41"};
{\ar@{-} "v12"; "v22"};
{\ar@{-} "v12"; "v32"};
{\ar@{-} "v12"; "v42"};
{\ar@{-} "v13"; "v23"};
{\ar@{-} "v13"; "v33"};
{\ar@{-} "v13"; "v43"};
{\ar@{-}@[green]@[thicker] "v14"; "v24"};
{\ar@{-}@[green]@[thicker] "v14"; "v34"};
{\ar@{-}@[green]@[thicker] "v14"; "v44"};
{\ar@{-} "v15"; "v25"};
{\ar@{-}@[green]@[thicker] "v15"; "v35"};
{\ar@{-}@[green]@[thicker] "v15"; "v45"};
{\ar@{-} "v16"; "v26"};
{\ar@{-} "v16"; "v36"};
{\ar@{-} "v16"; "v46"};
{\ar@{-} "v17"; "v27"};
{\ar@{-} "v17"; "v37"};
{\ar@{-} "v17"; "v47"};
{\ar@{-}@[green]@[thicker] "v24"; "v34"};
{\ar@{-}@[green]@[thicker] "v24"; "v44"};
{\ar@{-}@[green]@[thicker] "v35"; "v45"};
{\ar@{-} "v46"; "v37"};
"v47"; "v26" **\crv{ (-7, 130) & ( -30, 93) };
"v47"; "v33" **\crv{ (25, 103) & (30, 85) & (25, 67) };
"v27"; "v45" **[green][thicker]\crv{ (-20, 80) & (-7, 80) };
"v27"; "v35" **[green][thicker]\crv{ (7, 103) & (15, 90) & (18, 67) & (9, 67) };
"v37"; "v26" **\crv{ (-7, 108) & (-7, 90) };
"v46"; "v25" **\crv{ (-20, 98) & (-25, 87) };
"v36"; "v25" **\crv{ (-7, 90) };
"v36"; "v42" **\crv{ (22, 62) & ( 25, 44) & (32, 26) & (10, 24) };
"v44"; "v31" **\crv{ (15, 70) & (25, 36) & (25, -5) & (7, -5) };
"v34"; "v41" **\crv{ (3, 31) & (6, 18)};
"v43"; "v22" **\crv{ (-20, 40) & (-25, 32) & (-15, 20) };
"v43"; "v32" **\crv{ (20, 40) & (25, 32) & (15, 20) };
"v23"; "v41" **\crv{ (-30, 5) };
"v23"; "v31" **\crv{ (-7, 22) & (30, 12) };
"v33"; "v22" **\crv{ (-7, 31) };
"v42"; "v21" **\crv{ (-9, 0) };
"v32"; "v21" **\crv{ (-6, 14) };
};
\endxy\]
\caption{Two cospectral $3$-regular graphs constructed by Stark and Terras in \cite{ST2} by zeta function and covering considerations.}
\label{fig:TerrasCospectral}
\end{figure}

Figure \ref{fig:TerrasCospectral} is an example of two graphs $X_1$ and $X_2$ which Stark and Terras constructed as a consequence of zeta function and covering considerations in \cite{ST2}.  These graphs were constructed to have the same Ihara-Selberg zeta function and are thus cospectral as well since they are $3$-regular.  It is fairly straight forward to check that they are not isomorphic; however, we will use our zeta function to prove this, showing how we can get more leverage by controlling paths more precisely.

Both $X_1$ and $X_2$ have exactly 4 triangles.  We see this explicitly by considering the coefficients of their characteristic polynomials as in \cite{Biggs} or by noting that the coefficient of $u^3$ in $Z_{X_1}^{-1}(u)$ is $-8$, which is minus twice the number of triangles in $X_1$ \cite{Storm}.  We can find the triangles quickly by inspection.  In $X_1$, we've singled out two disjoint triangles in red.  In $X_2$, all four triangles are represented in green.

We now suppose that $X_1$ and $X_2$ are isomorphic.  We change $X_1$ into a hypergraph by replacing the red triangles with hyperedges on their vertices.  We can now compute the zeta function for this hypergraph.  As before, we've restricted the prime cycles on $X_1$ by throwing out any prime cycle that uses two red edges in a row.

If $X_1$ and $X_2$ are isomorphic, we should be able to repeat the transformation from graph to hypergraph in $X_2$ and finish with isomorphic hypergraphs.  There are four possible ways to create a hypergraph from $X_2$ in the same manner as we did for $X_1$.  For each green subgraph, we have a choice of two triangles to focus on, and there are two such green subgraphs.

Now a simple comparison of generalized Ihara-Selberg zeta functions distinguishes the graphs.  All four of the hypergraphs constructed from $X_2$ actually have the same zeta function.  However, the hypergraph we constructed from $X_1$ has a different zeta function.  Hence, these two graphs are not isomorphic.  Thus by making our paths more specific to the structure of the graphs, we've actually managed to get around Quenell's result and give a zeta function proof that $X_1$ and $X_2$ are not isomorphic.

We should mention that there is a drawback to this method as well.  We were fortunate that our example had a relatively small number of triangles.  As the number of non-disjoint triangles grows, we have to consider more and more potential hypergraphs.  Here, we only had to consider $4$ potential hypergraphs constructed from $X_2$; however, this was a graph with quite a small number of triangles.  Other options would be to make every possible triangle into a hyperedge; then, you would only have to compare one generalized Ihara-Selberg zeta function for each initial graph.  For this example, changing all four triangles into hyperedges and then computing the generalized Ihara-Selberg zeta function of the resulting hypergraphs also distinguishes the graphs.  We hope to explore these methods more at a later date.

All computations of zeta functions referenced in this section are available from the author by request.

\newpage
\bibliographystyle{plain}
\bibliography{../MainBiblio}
\end{document}